\documentclass[twoside,12pt]{article}
\def\ifundefined#1{\expandafter\ifx\csname#1\endcsname\relax}
\makeatletter \usepackage{theorem} {\AtBeginDocument{ \makeatletter
    \pagestyle{myheadings} \markboth{\hfill\ifundefined{@authorshort}
      \@author \else \@authorshort \fi\hfill}
    {\hfill\ifundefined{@titleshort} \@title \else \@titleshort
      \fi\hfill} \makeatother } }

\ifundefined{thm} \theorembodyfont{\slshape}
\newtheorem{thm}{Theorem}[section] 
\newtheorem{lem}[thm]{Lemma} \newtheorem{cor}[thm]{Corollary}

\theorembodyfont{\upshape} 

\newtheorem{example}[thm]{Example}

\newtheorem{rem}[thm]{\mdseries\scshape Remark} \fi

\ifundefined{proof}

\newenvironment{proof}[1][\proofname]{\par \normalfont
  \topsep6\p@\@plus6\p@ \trivlist
\item[\hskip\labelsep\scshape #1{.}]\ignorespaces
}{%
$\qed$\endtrivlist } \newcommand{\proofname}{Proof} \fi

\let\@ldmaketitle=\maketitle \renewcommand{\maketitle}{{\def\newpage{}
    {\footnotesize \parbox[t]{0.45\textwidth}{Preprint
    \texttt{math/9902012}, 1999}
      \hfill \parbox[t]{0.45\textwidth}{\whereappear}}
    \@ldmaketitle}} 
\newcommand{\whereappear}{
  Odessa University Herald, \textbf{3}(1998), pp.~60--63.}

\providecommand{\dedicatory}[1]{}
\providecommand{\keywords}[1]{\begingroup \def \protect {\noexpand
    \protect \noexpand }\xdef \@thefnmark { }\endgroup
  \@footnotetext{{\em Keywords and phrases.\/} #1}}
\providecommand{\AMSMSC}[2]{\begingroup \def \protect {\noexpand
    \protect \noexpand }\xdef \@thefnmark { }\endgroup
  \@footnotetext{{1991 \it Mathematical Subject Classification.\/}
    Primary: #1; Secondary: #2.}}

\newcommand{\authorshort}[1]{\gdef\@authorshort{#1}}
\newcommand{\titleshort}[1]{\gdef\@titleshort{#1}}

{\@definecounter{equation}}{\@newctr {equation}[section]}

\def\p@enumi{\thethm.}
\def\ifundefined#1{\expandafter\ifx\csname#1\endcsname\relax}
 \newcommand{\comment}[1]{}
 \usepackage{amsfonts}

\newcommand{\algebra}[1]{\ensuremath{\mathfrak{#1}}}
\newcommand{\Heisen}[1]{\ensuremath{{\mathbb{H}^{#1}}}}
\newcommand{\hHeisen}[1]{\ensuremath{{{ \widehat{\mathbb H}^{#1}}}}}

\newcommand{\object}[2][\,]{\ensuremath{\mathrm{#2}#1}}
\newcommand{\Space}[2]{\ensuremath{ {\mathbb{#1}^{#2}} }}
\newcommand{\such}{\,\mid\,}

\newcommand{\FSpace}[2]{{\ensuremath{ #1_{#2} }}}

\newcommand{\bos}{\mathcal{B}} \ifundefined{qed}
\DeclareMathSymbol{\qed}{0}{AMSa}{"03} \fi
\newcommand{\fourier}[1]{\ensuremath{\mathcal{F}^{#1}}}

 \newcommand{\Cstar}{$C^{*}$}

\newcommand{\norm}[1]{\left\| #1 \right\|}
\newcommand{\modulus}[1]{\left|#1\right|}
\newcommand{\scalar}[2]{\langle #1,#2\rangle}

\providecommand{\eqref}[1]{\textup{(\ref{#1})}}

\renewcommand{\baselinestretch}{1.3}
\newcommand{\rref}[1]{\textup{\ref{#1}}}
\newcommand{\Hspace}{\mathcal{ H}}
\makeatother
\title{Harmonic Analysis and Localization Technique}

\author{ Vladimir V. Kisil} \date{}
\begin{document}
\pagestyle{myheadings} \markboth{\hfill V. Kisil \hfill} {\hfill
  Harmonic Analysis and Localization Technique\hfill}

\maketitle \comment{\newpage
  \emph{Abbreviated title}: Harmonic and Local Analysis\\
  \newpage }
\begin{abstract}
  We study relationships between different formulations of the
  \emph{local principle}. Also we establish a connection among the
  local principle and the non-commutative Fourier transform approach
  to the investigation of convolution operator algebras.
  
  \keywords{Sectional Representations, Local Principle, Fourier
    Transform, Convolution Algebras.}
  
  \AMSMSC{46L45, 22E25}{47C15, 47D50}
\end{abstract}
\section{Introduction}

This paper is devoted to the {\em local principle\/}, which is a
canonical method for the study of operator algebras.
Theorem~\ref{th:relative} establish the coincidence of the easy to use
local principle based on the existence of a central commutative
subalgebra~\cite{Douglas72} and the general local principle
constructed on a set of ideals~\cite{Hofmann72}.

The second question of this paper is the correspondence between the
local principle and the non-commutative Fourier transform.  If we have
a group $G$ with the group operation $*$ and the Haar ($=$ invariant)
measure $dg$, it seems quite natural to introduce a group algebra
\algebra{G} associated with the group $G$.  Now {\em the
  noncommutative Fourier transformation\/} established a direct
connection between representations of the group algebra \algebra{G}
and representations of the group $G$ itself (see
\cite{Kirillov76,MTaylor86}).

On the other hand, the convolution algebra \algebra{G} has a regular
representation in the space of bounded operators
$\bos(\FSpace{L}{2}(G))$. This representation is introduced as the
integral of the shift operator $\pi_r(g)$
\begin{displaymath}
[\pi_r(g)f](h)=f(h*g)
\end{displaymath}
giving rise to the regular representation of the group $G $ on the
space $\FSpace{L}{2}(G)$:
\begin{equation}
\pi_r:\algebra{G}\rightarrow \bos(\FSpace{L}{2}(G)):  k(g)\mapsto
K=\int_{G} k(g)\pi_r(g)\,dg,
\label{eq:rconv}
\end{equation}
Therefore {\em the local principle\/} from the operator theory can be
applied to the operator algebra \algebra{G}.  The two ways lead to the
same answer, so a relation between the noncommutative Fourier
transformation and the local principle should exist. Note, that
harmonic analysis is described in terms of the Plancherel
\emph{measure} and sectional representations are constructed in field
$*$-\emph{topology}.  The correspondence between these two different
mathematical objects is clarified in Theorem~\ref{th:correspond}. An
application of the achieved result to the support of the Plancherel
measure is given in Corollary~\ref{co:center}. All proofs are easy and
they are skipped.

\comment{
\section{Preliminaries}
We introduce the main definitions and notations for future use.

\subsection{The Noncommutative Fourier Transformation}
\label{ss:Fourier}

There are many excellent books on this subject (see for example,
\cite{Kirillov76,MTaylor86}), thus we only fix some notations for
completeness here. Let $G $ be an unimodular type 1 group and
$\widehat{G}$ be the set of the all (up to unitary equivalence)
irreducible unitary representations $\pi_i$ of the group $G$ on
Hilbert spaces ${\Hspace_{i}}$.  For arbitrary $\pi_i \in \widehat{G}$
and $k\in \FSpace{L}{1}(G) $ the following operator on $\Hspace_{i} $
is defined (compare with \eqref{eq:rconv}):
\begin{equation}\label{eq:pigroup}
\pi_i (k) =\int_{G}k(g)\pi_i(g) dg,
\end{equation}
The usual way to define such operator and to proof its boundedness is
the following: for arbitrary $\eta, \nu \in \Hspace_{i}$ and
$k(g)\in\FSpace{L}{1}(G)$ the integral
\begin{displaymath}
\scalar{\eta}{\pi_i (k)\nu} =\int_{G}k(g)\scalar{\eta}{\pi_i(g)\nu} dg
\end{displaymath}
is obviously well defined (here $\scalar{\cdot}{\cdot}$ is the scalar
product in $\Hspace_{i} $) and so in turn it defines an operator
$\pi_i(k) $ from $\bos(\Hspace_{i})$.  A simple calculation shows
again that the identity \eqref{eq:pigroup} defines a unitary
representation of the group algebra \algebra{G} with
$\FSpace{L}{1}(G)$ kernel for the group $G$:
\begin{equation}\label{eq:pialgebra}
\pi_i:\ (\hbox{convolution with kernel }k\in \FSpace{L}{1}(G)) \mapsto
\pi_i (k).
\end{equation}
\noindent Note also that, if $G$ is an unimodular type 1 group (for
example a nilpotent Lie group), then there exists the Plancherel
measure $d\mu(\pi)$ on $\widehat{G}$ and the Plancherel identity
holds:
\begin{equation}\label{eq:Plancherel}
\norm{k} ^{2}_{L^{2}(G)}=\int_{\widehat{G}}\norm{\pi_i(k)}^{2}_{HS}
d\mu
(\pi ).
\end{equation}
Here $\norm{\pi(k)}^{2}_{HS}=\object{tr}(\pi_i(k)^*\pi_i(k)) $ is the
squared Hilbert--Schmidt norm of $\pi_i(k)$. The Plancherel identity
provides a ground for the reconstruction of the group algebra
\algebra{G} from representations $\pi_i(\algebra{G}) $.  It is worth
to mention here that {\em the support\/} of the Plancherel measure may
be quite different from $\widehat{G}$. The following Example
illustrates the common case (see also Corollary~\ref{co:center}).
\begin{example}\label{ex:heisen} Let us consider the Heisenberg group 
  \Heisen{n}~\cite{Kirillov76,MTaylor84,MTaylor86}. All irreducible
  representations of the group \Heisen{n} are given (up to unitary
  equivalence) by the Stone-von Neumann theorem. For any $\lambda \in
  (0, \infty)$ the irreducible unitary representations on
  $\FSpace{L}{2}(\Space{R}{n})$ are given by
\begin{equation}\label{eq:heisen1}
\pi_{\pm\lambda}(t, x, y)=e^{i(\pm\lambda tI\pm\lambda^{1/2}
yX+\lambda^{1/2}xD)},
\end{equation}
where $yX$ and $xD$ are:
\begin{eqnarray*}\displaystyle
(yX)u(z)&=&\sum_{j=1}^n y_{j}z_{j}u(z), \\
(xD)u(z)&=&(\frac{1}{i}) \sum_{j=1}^n x_{j}\frac{\partial u}{\partial
z_{j}}.
\end{eqnarray*}
For $(p, q)\in\Space{R}{2n}$, there are also one-dimensional
representations
\begin{equation}\label{eq:heisen2}
\pi_{(p, q)}(t, x, y) u=e^{i(px+qy)}u, \ u\in\Space{C}{}.
\end{equation}
So $\hHeisen{n}$ may be identified with $(\Space{R}{}\setminus 0 )\cup
\Space{R}{2n}$. But somewhat surprising properties of the Plancherel
measure are the following: $d\mu(\Space{R}{2n})=0,\ 
d\mu(\Space{R}{}\setminus 0)= \modulus{t}^n dt$. One can find an
explanation in Remark~\ref{re:nilpotent}.
\end{example}

\subsection{Basics of the Sectional Representation Theory} \label{ss:local}

Sectional representations theory is attempt to construct a
Gelfand-like isomorphism for non-commutative algebras. The detailed
account can be found in~\cite{DaunHof68,Hofmann72}, we present only
main ideas in connection with the theory of \Cstar--algebras.

To describe a topological algebra \algebra{R} one can take a family
$f_b:\algebra{R}\rightarrow\algebra{R}_b$, $b\in B$ of surjective
homomorphisms. Let $E$ be a disjoint union of $\algebra{R}_b$. The
crucial point is such topologies (called field $*$-topology) on $E$
and $B$ that
\begin{enumerate}
\item The projection $\pi: E \rightarrow B$, such that $\pi(x)=b$ for
  every $x\in \algebra{R}_b$, is continuous.
\item All global sections $\widehat{a}: B \rightarrow E$, such that
  $\widehat{a}(b)=f_b(a)$ are continuous.
\item The topologies induced on $R_b$ by $R$ (through $f_b$) and $E$
  coincide.
\end{enumerate}
In such topology image of $\algebra{R}$ belong to the algebra of
continuous bounded global section $\Gamma^b(B)$ over $B$. Note, that
one can alternatively speak not about homomorphisms $f_b$ but about
primitive ideals ${\cal I}_b$ corresponding to them.

For the \Cstar-algebras \algebra{R} one can start from an arbitrary
subset $B\in \object{Prim}\algebra{R}$ of set of primitive ideals and
produce the canonical sheaf of algebras. To obtain injection
$\algebra{R}\rightarrow\Gamma^b(B)$ the obvious condition $\bigcap
B=0$ should be imposed. The finest possible decomposition is obtained
by
\begin{thm}\textup{\cite[8.6]{Hofmann72}} Let \algebra{R} be a
  \Cstar-algebra with identity and let $\object{Prim'}\algebra{R}$ the
  space of its primitive ideals $\object{Prim}\algebra{R}$ with the
  coarsest topology making the function
  $\object{Prim}\algebra{R}\rightarrow \object{Spec} {\cal Z}:
  I\mapsto I\cup {\cal Z}$ with the centrum $Z$ of \algebra{R}
  continuous. Then there is a field of \Cstar-algebras over
  $\object{Prim'}\algebra{R}$ with stalks $\algebra{R}/I'$, $I'=\{a\in
  \algebra{R}\such \lim\norm{a+J}=0$ as J approaches $I$ in
  $\object{Prim'}\algebra{R}\}$ such that the Gelfand isomorphism of
  $\algebra{R}$ into \Cstar-algebra of global sections in this field
  is an isomorphism.
\end{thm}
\begin{rem}\label{re:regular}
  The injection condition ($\bigcap B=0$) is satisfied for any subset
  $B$ dense in $\object{Prim'}\algebra{R}$. Then
  $\object{Prim'}\algebra{R}$ is the unique (up to natural
  equivalence) Stone-\v Cech compactification (Hausdorffization) of
  $B$. Particularly it is rich enough to separate any point from a
  closed set not containing the point by a continuous function into
  $[0,1]$~\cite[8.4]{Hofmann72}.
\end{rem}
}

\section{A System of Ideals and the Algebra Center}\label{se:center}

To apply the local principle from \cite{Hofmann72} one need to find an
appropriate system of ideals. If the algebra under consideration has a
non-trivial abelian subalgebra then the answer is given by the
following

\begin{thm}\label{th:Douglas}\textup{\cite[Proposition
    4.5]{Douglas72}} If \algebra{R} is a \Cstar-algebra, ${\cal Z}$ is
  an abelian \Cstar--subalgebra contained in \algebra{R} having
  maximal ideal space $M_{\cal Z} $, and for $x $ in $M_{\cal Z} $,
  ${\cal I}_x $ is a closed ideal in \algebra{R} generated by the
  maximal ideal $\{Z\in{\cal Z}:\widehat{Z}(m)=0\}$, then
  $\vspace{5pt}\bigcap_{x\in M_{\cal Z}}\, {\cal I}_x=\{0\} $. In
  particular, if $\Phi_x$ is *-homomorphism from \algebra{R} to
  $\algebra{R}/{\cal I}_x$ then $\vspace{5pt}\sum_{x\in M_{\cal
      Z}}\bigoplus \Phi_x$ is *-isomorphism of \algebra{R} into
  $\vspace{5pt}\sum_{x\in M_{\cal Z}}\bigoplus \algebra{R}/{\cal
    I}_x$.  Moreover, $T$ is invertible in \algebra{R} if and only if
  $\Phi_x(T)$ is invertible in $\algebra{R}/{\cal I}_x $ for $x $ in
  $M_{\cal Z}$.
\end{thm}

In the case of a trivial (or unsuitable) center of the algebra
\algebra{R} the following trick is commonly used: one can introduce a
larger algebra $\widetilde{\algebra{R}}\supset\algebra{R}$ with
appropriate center $\widetilde{\cal Z}$ and use
Theorem~\rref{th:Douglas} for the algebra $\widetilde{\algebra{R}}$.
Then the description of the algebra \algebra{R} can be obtained as a
subalgebra of the algebra $\vspace{5pt}\sum_{x\in M_{\widetilde{\cal
      Z}}}\bigoplus \widetilde{\algebra{R}}/{\cal I}_x$.
\begin{example}\label{ex:bisingular} The algebra
  \algebra{B} of bisingular operators \cite{Simonenko74} has only
  operators of multiplication by a constant in its center.
  Nevertheless, one can introduce the largest algebra \algebra{R} of
  operators having the form $a(x)\fourier{-1}b(\xi)\fourier{}$ with
  $a(x)\in\dot{\Space{R}{2}}$ and $b(\xi)$ be a homogeneous degree
  zero function whose restriction on \Space{S}{1} is piecewise
  continuous.  This algebra contains the bisingular operators and have
  a central commutative (up to compact operators) subalgebra of
  Calderon-Mikhlin-Zigmund operators with continuous symbols. The
  application of Douglas' technique allows to obtain the full
  description of the algebra of symbols of the algebra \algebra{R}
  (and the algebra \algebra{B} correspondingly).
\end{example}

It is remarkable that this trick may be applied to the general local
principle. The following Lemma settles the particular case of regular
base $B$ in the field $*$-bundle topology.
\begin{lem}\label{le:relative}
  Let a \Cstar-algebra \algebra{R} have a system $B$ of primitive
  ideals with property $\bigcap B=0$ and $B$ be compact in the field
  $*$-topology.  Then there is a \Cstar-algebra
  $\widetilde{\algebra{R}}$ with non-trivial center $\widetilde{{\cal
      Z}}$ and an embedding of the algebra $\algebra{R}$ on
  $\widetilde{\algebra{R}}$ such that the application of the general
  local principle by the systems $B$ and the method from
  Theorem~\rref{th:Douglas} to the subalgebra $\widetilde{{\cal Z}}$
  give the same result.
\end{lem}
\comment{
\begin{proof} We  give the direct construction of the required algebra
  $\widetilde{\algebra{R}}$ based on the canonical \Cstar-bundle
  defined by the algebra \algebra{R} and by the system of ideals $B$.
  Let \algebra{C} be the algebra of sections having the form
  $a(t)e_t$, where $a(t)$ is a continuous (and thus bounded) function
  on $B$ and $e_t$ is identity from the algebra ${\cal R}_t$. Then the
  algebra $\widetilde{\algebra{R}}$ can be defined as the closed
  algebra generated by elements from \algebra{R} and \algebra{C}.
  Clearly \algebra{C} is a central commutative subalgebra of
  $\widetilde{\algebra{R}}$ having the set $B$ of its maximal ideals
  compact. It is not hard to check then that the Douglas' techniques
  and the general local principle coincide.
\end{proof}
}
\begin{rem}{}\label{re:bisingular-center} We come back to the
  bisingular operator algebra from Example~\ref{ex:bisingular} again.
  It is easy to verify it this case that the largest algebra of
  sections cannot be realized as an operator algebra on the initial
  Hilbert space $\FSpace{L}{2}(\Space{R}{2})$.
\end{rem}
\begin{thm}\label{th:relative}
  Let a \Cstar-algebra \algebra{R} have a system $B$ of ideals such
  that $\bigcap B=0$.  Then there is a \Cstar-algebra
  $\widetilde{\algebra{R}}$ with non-trivial center $\widetilde{{\cal
      Z}}$ and an embedding of the algebra $\algebra{R}$ on
  $\widetilde{\algebra{R}}$ such that the application of the canonical
  shelf construction by the systems $B$ and the method from
  Theorem~\rref{th:Douglas} to the subalgebra $\widetilde{{\cal Z}}$
  give the same result.
\end{thm}
\comment{
\begin{proof}
  As we know (see Remark~\ref{re:regular}) $B$ has a regularization
  $\object{Prim'}\algebra{R}$ points of which are again primitive
  ideals of \algebra{R}. Thus we can construct the canonical sectional
  representation associated with $\object{Prim'}\algebra{R}$ and apply
  the previous Lemma.
\end{proof}
}

\section{The General Local Principle and the Support of the Plancherel
  Measure}

Operator algebras in analysis usually have rich groups of symmetries,
which are a symmetry group of differential equations define the
analyticity property. Thus they can be treated as
convolutions~\cite{Kisil95a}. We would like to apply the localization
technique in such an environment.  Let $B$ be a family of primitive
ideals in the group \Cstar--algebra \algebra{G} with the only
property:
\begin{equation}\label{eq:zero}
\bigcap B = 0
\end{equation}
\noindent (for example, all ideals generated by all maximal  ideals of
a central commutative subalgebra of our algebra, see
Section~\rref{se:center}). Then the following quotient mapping for any
$ J_b$
\begin{equation}\label{eq:local}
\pi _{t}:(\hbox{convolution with kernel }k)\mapsto
{\hbox{(convolution
with kernel $k$)}}/{\hbox{$ J_b$}}
\end{equation}
\noindent generates a representation of the group algebra. Condition
\eqref{eq:zero} ensures that all irreducible representations of the
group algebra may be obtained as (sub)rep\-re\-sen\-ta\-tions of
\eqref{eq:local}. However, the careful extracting of irreducible
representations may be a rather difficult problem~\cite{Kisil94a}.

The following Lemma plays a fundamental role in the establishment of a
correspondence between the local principle and the Fourier transform.
\begin{lem}\label{le:support} Let $G$ be an unimodular type 1 
  exponential group, let $\widetilde{G}$ be the support of the
  Plancherel measure $d\mu $ in $\widehat{G}$. Then the family of
  two-sided ideals
\begin{displaymath}
 J(\pi_i)=\{K\mid  K \hbox{ is a convolution with kernel  $k$ such,
that
} \pi_i (k)=0\},\ \pi_i \in \widetilde{G}
\end{displaymath}
\noindent satisfies the condition \eqref{eq:zero}.
\end{lem}
\comment{
\begin{proof}
  Direct application of the Plancherel identity \eqref{eq:Plancherel}
  and the simple observation that the integral from
  \eqref{eq:Plancherel} may be taken only on the support of the
  Plancherel measure imply the assertion.
\end{proof}
}
Lemma~\ref{le:support} suggests that one may use this family $J(\pi_i)
$ of ideals for localization according to the general local principle.
It obviously follows from the definition of the family $J(\pi_i) $ of
ideals that the rules non-commutative Fourier transform and local
technique define just the same representations of our group algebra.
Thus in the mentioned case the noncommutative Fourier transformation
and the local principle give, in fact, the same description of the
group convolution algebra. This may be summarized as follows
\begin{thm} \label{th:correspond}
  Under the assumptions and notations of Lemma \rref{le:support} we
  have: All representations of the group algebra \algebra{G} given by
  the Fourier transform
\begin{displaymath}
\pi_i:\ (\hbox{convolution with kernel }k\in \FSpace{L}{1}(G)) \mapsto
\pi_i (k)
\end{displaymath}
are contained (as subrepresentation, possibly) within the
representations
\begin{displaymath}
\pi _{i}:(\hbox{convolution with kernel }k)\mapsto {\hbox{(convolution
with kernel $k$)}}/{\hbox{$ J(\pi_i)$}},
\end{displaymath}
where
\begin{displaymath}
J(\pi_i)=\{K\mid  K \hbox{ is a convolution with kernel  $k$ such,
that
} \pi_i (k)=0\},\ \pi_i \in \widetilde{G}.
\end{displaymath}
The support of the Plancherel measure $\widetilde{G}$ is a dense
subset of $\object{Prim'}\algebra{C}$ in field $*$-topology.
\end{thm}
After the remarks made the proof is not necessary. The next result
shows how the described correspondence can be applied.
\begin{cor}\label{co:center}
  Under the assumption of Theorem \rref{th:correspond}, the Plancherel
  measure of $\widehat{G}$ is supported on those representations of
  $G$ which do not contain the center of $G$ in their kernels.
\end{cor}
The author was partially supported by the INTAS grant 93--0322.
During the final preparation of the paper the author enjoyed the
hospitality of Universiteit Gent, Varkgroep Wiskundige Analyse,
Belgium.  \small \renewcommand{\baselinestretch}{0.9}
\newcommand{\noopsort}[1]{} \newcommand{\printfirst}[2]{#1}
\newcommand{\singleletter}[1]{#1} \newcommand{\switchargs}[2]{#2#1}
\newcommand{\irm}{\textup{I}} \newcommand{\iirm}{\textup{II}}
\newcommand{\vrm}{\textup{V}}

\end{document}